\documentclass[english,12pt]{article}

\usepackage{graphicx}
\usepackage{amssymb,amsmath}
\newcommand{\D}{{\cal D}}
\newcommand{\swsw}{{\tt ss}}
\newcommand{\mama}{{\tt mm}}
\newcommand{\swma}{{\tt sm}}
\newcommand{\masw}{{\tt ms}}

\newcommand{\doors}{{
\normalfont\initfamily
\fontsize{12mm}{12mm}\selectfont  \raisebox{1ex}{\small The}doors
\normalfont\initfamily}}

\input Starburst.fd

\newcommand*\initfamily{\usefont{U}{Starburst}{xl}{n}}

\newenvironment{myproof}{\noindent {\it Proof} }{$\Box$ }

\begin{document}
\title{\doors}

\author {Sasha Gnedin\\{\tt A.V.Gnedin@uu.nl }}

\maketitle

\begin{abstract}
\noindent
We emphasize the dominance in the Monty Hall problem,  both in the classical scenario and its multi-door generalization.
This is used to show optimality of the class of always-switching strategies for nonuniform 
allocation of the prize and arbitrary door-revealing mechanism in the event of match.
\end{abstract}


\paragraph{\large To  switch or not to switch$\ldots$}
You have two tries to find a prize hidden behind one of three doors.
You are first asked  to choose a door but not open it yet. 
Then one of the unchosen doors with no prize behind it will be revealed, 
and you will be offered a second try. 
You win if the final choice falls on the door concealing the prize.
{\it Does it matter if you switch?}

This is the famous Monty Hall problem.
It is commonly assumed, often implicitly, that the prize is hidden ``uniformly at random'',
meaning that it is equally likely to be behind each of the doors.
There are two basic arguments showing that switching is better.
The arguments are among the best known pieces of elementary probability theory
(see e.g. the article \cite{Gillman} 
), 
but we need to sketch them here to show our point.

The {\it cases argument} goes as follows. Suppose the prize is behind door $\theta$.
If you choose door 1 and then never switch,
you win if $\theta=1$. If you choose door 1 and   always switch, you win when $\theta\in\{2,3\}$. The cases are mutually exclusive,  
therefore the probability to win with switching is  $2/3$.

A major interest to the problem stems from the observation
that many people find this solution counter-intuitive, as they feel that the odds 
in favor of switching  
are 1:1 when one of the doors is revealed.
Defining these 
{\it conditional}~ odds rigorously requires a further assumption 
on the random mechanism selecting  a door to reveal in the event of a {\it match}, when 
the initial choice falls on the door hiding the prize.
If the randomization is performed by tossing the same fair coin,
the conditional odds in favor of switching are 2:1 with certainty.
This is the second, {\it odds argument}.

If the randomization is done by tossing a biased coin, with bias allowed to depend
on door $\theta\in \{1,2,3\}$ at which a match has occured, the odds are still never unfavorable for switching. 
See Wikipedia pages on the Monty Hall problem as a general source of references,
\cite{Gill1} for discussion of assumptions, \cite{Gill2}  
for a  summary of other arguments and their interconnections 
and \cite{Rosenhouse} for many variations of the problem.


In the nonuniform case, denoting $p_\theta$ the probability to find the prize behind door $\theta$ ($\theta=1,2,3$),
and $q$ the conditional probability to leave  door 2 unrevealed in the event of match at door $1$,
the odds in favor of switching from door $1$ to $2$ are $p_2: p_1q$, hence switching 
is  advantageous
when $p_1q<p_2$. The  inequality holds  always for $p_1<p_2$ but depends on $q$ if $p_1\geq p_2$.   
For instance, 
if $p_1=4/9, ~p_2=3/9, ~p_3=2/9$ switching from door $1$ to $2$ is better exactly when $q<3/4$.
Despite disadvantage in some situations, 
switching cannot be devaluated: in the example
the initial choice of  door 1 was not optimal, and  you would get a higher chance to 
win the prize by  first choosing door 3 then switching  all the time as the second try is offered.
When the probabilities
$p_\theta$ are not equal it is useful to evaluate the switching action in combination with the initial choice.

\paragraph{\large $\ldots$ how can you know to which?}
Without assigning probabilities  to the doors the question
{\it Does it matter if you switch?} is sometimes regarded  as not well-posed mathematically (e.g. \cite{Ethier}, p. 56).
Nevertheless, the question makes sense 
if we consider the decision-making
as a two-step process and compare the {\it class} of strategies that always switch
with all other strategies. It turns that the problem has {\it dominance} that speaks
in favor of always-switching strategies.

To see the dominance we need to ``switch the door'' in the two basic arguments.
Let us start with the cases argument, and  compare two {\it constant-action} strategies 
$A=$``choose door 1, never switch'' with $B=$``choose door 2, always switch''. 
They might seem incomparable, because when $A$ is played door 1 cannot be revealed, and when $B$ is played door 2 cannot be revealed.
However, this is irrelevant and, obviously enough, $A$ wins for $\theta=1$ while $B$ wins for $\theta\in \{1,3\}$.
Thus $B$ is never worse than $A$ and if door $\theta=3$ is possible (in some sense), then $B$ is even better (in this very sense).

The dominance reasoning  parallel to the odds argument is slighly more complicated.
To increase generality suppose $\cal D$ is a finite set of doors with at least three elements.
You wish to guess door $\theta$ hiding the prize. You first choose door $x$ and then 
you are offered a switch to door $y\in\D\setminus\{x\} $. Both $x$ and $y$ can be $\theta$, all other doors have been revealed as useless. 
Now you need to decide whether $\theta$ and $x$ 
match or mismatch, that is whether $x=\theta$ or $x\neq \theta$.
A strategy is therefore a combination of  $x$ and a decision function $a(x,y)$ with values in the two-point set of actions $\{{\tt match, switch}\}$,
with {\tt match} meaning staying {\it with} $x$ and {\tt switch} meaning  switching {\it from} $x$.

We denote $a^*$ the always-switch decision function, with the only action $a^*(x,y)={\tt switch}$  for all $y\neq x$.

Your strategy is rewarded by means of the win-or-nothing payoff function $W$
defined as 
$$W(\theta,x,{\tt match})=1(x=\theta),~~~W(\theta,x,{\tt switch})=1(x\neq\theta),$$
where $1(\cdots)$ equals $1$ if $\cdots$ is true, and equals $0$ otherwise.
Two actions might seem permutable, but the symmetry is fallacious, because
for every $\theta$ there is only one $x=\theta$ and
at least two $x\neq \theta$.


The way we write the payoff function does not involve explicitly the dependence on the door
to which the switch is offered. It is this feature which will enable us to compare the always-switching strategies with others.
When switching is offered to some door $y\neq x$,  the inequality $x\neq\theta$ means, of course, $\theta=y$
as we know that other doors have no prize. 
Nevertheless, the intepretation  of the win with {\tt switch} as
a mismatch of the prize door with your initial choice 
is  more insightful and general. We may think, for example, of the game in which 
less doors are revealed as empty and
switching is offered 
to a subset of $\D\setminus\{x\}$, with the convention that the action is successful if
$\theta$ is there.

The following key lemma is obvious from the definitions, and the dominance is a consequence.

\vskip0.2cm
\noindent
{\bf Lemma} {\it ~For all  $\theta$ and  distinct $x, x'\in \D$
$$W(\theta,x,{\tt match})=1~~ \Rightarrow~~  W(\theta,x',{\tt switch})=1.$$
}
\vskip0.2cm
\noindent
{\bf Dominance Theorem}{\it ~ Suppose for some fixed $x\neq y'$ a decision function $a(x,\cdot)$ satisfies $a(x,y')={\tt match}$.
Then  
$$W(\theta,x,a(x,y))\leq W(\theta,y',{\tt switch})$$
for all $\theta$  and $y\neq x$, meaning that the strategy $(x,a(x,\cdot))$ is weakly dominated by $(y',a^*)$.}\\
\vskip0.2cm
\begin{myproof}
Strategy $(y',a^*)$ always wins unless $\theta=y'$, when $W(y', y',{\tt switch})=0$, but then 
$W(y',x,a(x,y'))= W(y',x,{\tt match})=  0$ as well.
\end{myproof}
\vskip0.2cm

One insightful way to explain the phenomenon of dominance is to observe that for every strategy $(x, a_x(\cdot))$ there exists
a `unlucky door' $u$ such that the strategy misses the prize when it is behind $u$, no matter in which admissible way
the switching option is offered. 
Then the strategy $(u, a^*)$ is weakly dominating $(x, a(x,\cdot))$.

We display next the payoff structure for all strategies
in the 3-door case, $\D=\{1,2,3\}$.
The notation must be self-explaining, but keep in mind that the value of variable $y$ 
in the upper row is given in the case of match $x=\theta$ (otherwise $y=\theta$).
For instance, let us check the entry 2\masw/2,1: the
notation  $2\masw$ encodes the strategy $x=2$, $a(2,3)={\tt switch}, a(2,1)={\tt match}$, so if $\theta=2, y=1$ this is a win,
$W(2,2,{\tt match})=1$. 
For entry  2\masw/1,3 we have the same decision function, $\theta=1$, in the notation of column 1,3 digit   `3' is  irrelevant
since we have a mismatch, so $y=\theta=1$ and $W(1,2,a(2,1))=W(1,2,{\tt match})=0$.

\begin{center}
\begin{tabular}{c|cccccc}
$\theta,y=$    & 1,2 & 1,3 & 2,1 & 2,2 & 3,1 & 3,2\\

\hline
1\swsw &0  & 0 &  1  &  1 & 1  &1\\
1\masw &1  & 0 &  0 &  0 & 1  &1\\
1\swma &0  & 1 &  1  &  1 & 0  &0\\
1\mama &1  & 1 &  0  &  0 & 0  &0\\
    &   &   &     &    &    &  \\
2\swsw &1  & 1 &  0  &  0 & 1  &1\\
2\masw &0  & 0 &  1 &  0 & 1  &1\\
2\swma &1  & 1 &  0  &  1 & 0  &0\\
2\mama &0  & 0 &  1  &  1 & 0  &0\\
    &   &   &     &    &    &  \\
3\swsw &1  & 1 &  1  &  1 & 0  &0\\
3\masw &0  & 0 &  1 &  1 & 1  &0\\
3\swma &1  & 1 &  0  &  0 & 0  &1\\
3\mama &0  & 0 &  0  &  0 & 1  &1\\
\end{tabular}
\end{center}
The dominance of always-switching strategies is seen by comparing the rows.

The payoff matrix appears in 
\cite{Haggstrom}. In \cite{Mondee} reduction by dominance was used to arrive at the  minimax winning probability $2/3$
(see e.g. \cite{Haggstrom},\cite{Gill1}).

Although the domination is {\it weak}, in the sense that dominated strategy cannot be {\it strictly} improved for all $\theta$,
it is  a serious ground to discard dominated strategies in many settings of decision making.
Think for example of guessing the right door when the prize is hidden by some algorithm.
You would not use a strategy if there is another one performing at least as good and in some situations even better.

\paragraph{\large To maximize your score $\ldots$}
In the Bayesian setting of the guessing problem probabilities are assigned to all 
values of the variables out of your control. One random variable is the door with the prize $\Theta$,
with some 
specified probabilities $p_\theta$ for each value $\theta\in \D$. Another random variable, the door offered for switching $Y$
is defined conditionally on the value $\Theta=\theta$ and your choice $x$.
In the event of mismatch $\Theta\neq x$ we have  $Y=\Theta$, and given $\Theta=x$ the variable $Y$ assumes
each admissible value $y\in \D\setminus\{x\}$ with some probability $q_{x,y}$.
Think of  biased roulette wheels, distinct for each door $x$ and each having $\#\D-1$ positions to generate $Y$.

The probability to win with strategy $a(x,\cdot)$ is equal to the expected payoff
\begin{eqnarray*}
{\mathbb E}W(\Theta,x,a(x,Y))&=&\\
\sum_{\theta\in\D\setminus \{x\}} p_\theta W(\theta,x,a(x,\theta))+\sum_{y\in\D\setminus\{x\}}p_x q_{x,y}W(x,x,a(x,y))&=&\\
\sum_{\theta\in\D\setminus \{x\}} p_\theta 1(a(x,\theta)={\tt switch})
+\sum_{y\in\D\setminus\{x\}}p_x q_{x,y}1(a(x,y)={\tt match})&=&\\
\sum_{y\in\D\setminus \{x\}} {\big\{} 
p_y 1(a(x,y)={\tt switch})
+p_x q_{x,y}1(a(x,y)={\tt match}){\big\}}&=&\\
\sum_{y\in\D\setminus \{x\}} {\big\{} 
p_y 1(a(x,y)={\tt switch})
+p_x q_{x,y}(1-1(a(x,y)={\tt switch}){\big\}}&=&\\
p_x+
\sum_{y\in\D\setminus \{x\}} {\big\{} 
(p_y -     p_x q_{x,y}) 1(a(x,y)={\tt switch}){\big\}}.
\end{eqnarray*}
In particular, the winning probability is 
 equal to $1-p_x$ for the always-switching $a^*$, and is equal to $p_x$ for always-matching $a(x,\cdot)\equiv{\tt match}$, 
as it was clear without computation.

\paragraph{\large $\ldots$ try first to miss \boldmath$the$ door!}
Finding the strategy with the highest winning chance (aka {\it Bayesian}~ strategy) 
seems cumbersome from the 
above explicit formula for ${\mathbb E}W(\Theta,x,a(x,Y))$. 
The odds in favor of switching from $x$ to $y$ are $p_{y}: p_x q_{x,y}$,
and these could be arbitrary,  so the sign of $p_y -     p_x q_{x,y}$ depends on $q_{x,y}$ 
if $p_x> p_y$.

The dominance helps. 
Discarding by dominance all strategies except $(x,a^*)$ we are left with maximizing $1-p_x$ over $x\in \D$.
Let $\theta^*$ be the least likely door to conceal the prize (such door need not be unique), 
$p_{\theta^*}=\min_{\theta\in\D}p_\theta$. 
The always-switching policy $(\theta^*,a^*)$ yields the highest winning probability $1-p_{\theta^*}$, that is
$$\max_{x\in\D,\,a(x,\cdot)}{\mathbb E}W(\Theta,x,a(x,Y))= {\mathbb P}(\Theta\neq\theta^*)=
1-p_{\theta^*}\,,$$
whichever the random rule to reveal the doors in the event of match.

Switching from  the least likely door is always beneficial since 
$$p_{y}\geq p_{\theta^*}\geq p_{\theta^*} q_{\theta^*,y}.$$ 
This could be concluded without evaluation of odds, 
directly from the optimality of the strategy $(\theta^*,a^*)$, 
which instructs that to get the prize it is optimal to first miss it with the highest possible probability.

\vskip0.4cm
\noindent

\noindent
{\small Instances of the problem with nonuniform distribution (3-door case) appear 
in \cite{Rosenhouse},\cite{Tijms}, \cite {Snell}.

Dominance was used in some multi-door game variations of the problem \cite{Bailey}. 
Apparently, the dominance (in which sense?) in the standard MHP 
 was conjectured  by columnist John Kay}
{\small
\begin{itemize}
\item[]{
 If my conjecture has been right switching would be a (weakly) dominant strategy, i.e. you can never lose by switching but you may gain}
({\it Financial Times}, 31 August 2005). 
\end{itemize}
}
\vskip0.2cm
\noindent
{\bf Acknowledgments} The author is indebted to Richard Gill and Jim Pitman whose polar views on the author's engagement in 
the MHP generated this note.

\end{document}